\documentclass[letterpaper, 10 pt, conference]{ieeeconf}  

\IEEEoverridecommandlockouts                              

\usepackage{cite}
\usepackage{color}
\usepackage[pdftex]{graphicx}
\graphicspath{{fig/}{jpeg/}}
\usepackage[cmex10]{amsmath}
\usepackage{amssymb}
\usepackage{algorithm}
\usepackage{algorithmic}
\usepackage[font=small]{caption}
\usepackage{subcaption}
\usepackage{setspace}
\usepackage{hyperref}
\usepackage{comment}
\usepackage{tikz}

\usepackage{needspace}

\usepackage{theorem}
\newtheorem{theorem}{Theorem}

\newtheorem{lemma}{Lemma}
\theoremstyle{definition}
\newtheorem{assumption}{Assumption}
\newtheorem{remark}{Remark}
\newtheorem{corollary}{Corollary}

\title{\LARGE \bf 
A Decentralized Quasi-Newton Method for Dual Formulations of Consensus Optimization
}

\author{Mark Eisen, Aryan Mokhtari, and Alejandro Ribeiro
\thanks{Work supported by NSF CAREER CCF-0952867 and ONR N00014-12-1-0997. Authors are with the Department of Electrical and Systems Engineering at the University of Pennsylvania, Philadelphia, PA 19104 USA. (email: {\tt maeisen, aryanm, aribeiro@seas.upenn.edu}).}
} %

\input{mysymbol.sty}

\begin{document} \maketitle

%
\begin{abstract}
This paper considers consensus optimization problems where each node of a network has access to a different summand of an aggregate cost function. Nodes try to minimize the aggregate cost function, while they exchange information only with their neighbors. We modify the dual decomposition method to incorporate a curvature correction inspired by the Broyden-Fletcher-Goldfarb-Shanno (BFGS) quasi-Newton method. The resulting dual D-BFGS method is a fully decentralized algorithm in which nodes approximate curvature information of themselves and their neighbors through the satisfaction of a secant condition. Dual D-BFGS is of interest in consensus optimization problems that are not well conditioned, making first order decentralized methods ineffective, and in which second order information is not readily available, making decentralized second order methods infeasible. Asynchronous implementation is discussed and convergence of D-BFGS is established formally for both synchronous and asynchronous implementations. Performance advantages relative to alternative decentralized algorithms are shown numerically.
\end{abstract}


%
\section{Introduction} \label{sec_intro}

We study the problem of decentralized consensus optimization where nodes of a network maximize a global objective function, while each of them has access to a different summand of the global objective function. To be more precise, consider a variable $\tbx\in\reals^p$ and a local strongly concave function $f_i: \mathbb{R}^p \rightarrow \mathbb{R}$ associated with node $i$. The goal of nodes is to solve the optimization problem 
\begin{equation}\label{original_problem}
\tbx^*:=\argmax_{\tbx\in\reals^p} f(\tbx) = \argmax_{\tbx\in\reals^p}\sum_{i=1}^n f_i(\tbx),
\end{equation}
while being allowed to exchange information with neighbors only. These problems arise in decentralized control
\cite{olfati2004consensus,Bullo2009,Cao2013-TII,LopesEtal8}, sensor networks
\cite{Schizas2008-1,KhanEtal10,cRabbatNowak04}, and
machine learning
\cite{bekkerman2011scaling,Tsianos2012-allerton-consensus,Cevher2014}.  

The theory and practice of first order methods to solve \eqref{original_problem} is well developed. There are multiple methods that solve \eqref{original_problem} in the primal domain \cite{nedic2009,nedic2010constrained,YuanQing,Shi2014} and a larger number of methods that solve \eqref{original_problem} through duality theory\cite{chatzipanagiotis2013augmented,rockafellar1976augmented,jakovetic2015linear,cRabbatNowak04,Jakovetic2014-1,Schizas2008-1}. However, and as is the case in centralized optimization, these first order methods are slow to converge when the objective function is ill-conditioned. This has motivated the development of decentralized second order methods which perform better than their first order counterparts, when the problems are not well conditioned and Hessians are available at reasonable computational cost \cite{NN-part1,bajovic2015newton,mokhtari2016decentralized}.
 Alas, evaluation and inversion of Hessians is a task that can be computationally impractical in some problems. When this is the case in centralized optimization, the solution comes in the form of resorting to quasi-Newton methods \cite{Broyden, Byrd, DingNocedal}. The goal of this paper is to develop a decentralized quasi-Newton method to handle problems that are not well conditioned and in which second order information is not readily available. 

We start by equating the solution of \eqref{original_problem} to the minimization of a suitable dual function (Section \ref{sec_problem_formulation}). A brief description of a regularized version of the centralized Broyden-Fletcher-Goldfarb-Shanno (BFGS) quasi-Newton method is then introduced (Section \ref{sec_bfgs_cent}). BFGS, regularized or not, can't be implemented in a decentralized manner because it relies on multiplying gradients by a curvature matrix that is not sparse. This limitation is overcome by the Decentralized (D-)BFGS method which relies on the observation that the appealing convergence traits of BFGS come from the curvature matrix satisfying a secant condition that can be expressed and satisfied in a decentralized manner (Section \ref{sec_dbfgs}). D-BFGS is a modification of regularized BFGS that maintains validity of this secant condition while ensuring the curvature matrix has a sparsity pattern matching the sparsity pattern of the graph. Asynchronous implementation of D-BFGS is further discussed (Section \ref{sec_async_dbfgs}). Convergence of D-BFGS is established for synchronous and asynchronous implementations (Section \ref{sec_convergence}) and performance advantages relative to alternatives are evaluated numerically (Section \ref{sec_numerical_results}). We close the paper with concluding remarks (Section \ref{sec_conclusion}). 

%
\section{Problem Formulation} \label{sec_problem_formulation}

We consider a decentralized system with $n$ nodes which are connected as per the graph $\ccalG=(\ccalV,\ccalE)$ where $\ccalV=\{1,\dots,n\}$ is the set of nodes and $\ccalE=\{(i,j)\}$ is the set of $m$ edges. We assume the graph $\ccalG$ is symmetric, i.e., $(i,j)\in\ccalE$ implies $(j,i)\in \ccalE $, and the graph does not contain self-loops. Further, define the neighborhood of node $i$ as the set $ n_i:=\{j\ |\ (i,j)\in\ccalE\}$ of nodes $j$ that are adjacent to $i$. The nodes have access to their local functions $f_i$ only and their goal is to find the optimal argument $\bbx^*  \in \mathbb{R}^{p}$ that minimizes the aggregate cost function, $\sum_{i=1}^n f_i(\tbx)$ in \eqref{original_problem}. To rewrite this optimization problem in a manner that is more suitable for decentralized settings, we introduce the variable $\bbx_i$ as a copy of the decision variable $\tbx$ kept at node $i$. We then rewrite the optimization problem in \eqref{original_problem} as
\begin{alignat}{2}\label{eq_primal_problem}
    \bbx^* \ :=\ &  \argmax_{\bbx_1,\dots,\bbx_n} \
                 && \sum_{i=1}^n f_i(\bbx_i) \nonumber\\
                 &  \ \text{s.t.}
                 && \bbx_i =  \bbx_j, \quad \forall\ (i,j) \in \ccalE.
\end{alignat}
Since the graph is connected, any feasible solution of \eqref{eq_primal_problem} satisfies $\bbx_1=\dots=\bbx_n$. With this restriction on the feasible set the cost functions in \eqref{original_problem} and \eqref{eq_primal_problem} become equivalent and the optimal argument $\bbx^*=\{\bbx_1^*,\dots,\bbx_n^*\}$ of \eqref{eq_primal_problem} has the form $\bbx_1^*=\dots=\bbx_n^*=\tbx^*$.

We tackle the solution of \eqref{eq_primal_problem} in the dual domain. Define then the dual variable $\bblambda_{ij}\in \reals^{p}$ associated with the constraint $\bbx_i=\bbx_j$ which is kept at node $i$. Moreover, define $\bblambda_i$ as the concatenation of all $\bblambda_{ij}$ for $j \in n_i$. Further, consider $\bblambda:=[\bblambda_1;\dots;\bblambda_n]\in\reals^{mp}$ as the concatenation of the $n$ dual variables $\bblambda_{i}$, to write the Lagrangian $\ccalL(\bbx,\bblambda)$ of the optimization problem in \eqref{eq_primal_problem} as 
\begin{equation}\label{eq_lagrangian}
\ccalL(\bbx,\bblambda) = \sum_{i=1}^n f_i(\bbx_i) + \sum_{(i,j) \in \ccalE} \bblambda_{ij}^T(\bbx_i - \bbx_j).
\end{equation}
Notice that the Lagrangian $\ccalL(\bbx,\bblambda) $ for a given dual vector $\bblambda$ is separable over the nodes. Hence, each node $i$ computes the local Lagrangian maximizer $\bbx_i(\bblambda)$ by solving the program 
\begin{equation}\label{eq_lagrangian_maximizers}
\bbx_i(\bblambda)=\argmax_{\bbx_i\in\reals^p} f_i(\bbx_i)+\sum_{j\in  n_i}(\bblambda_{ij}-\bblambda_{ji})^T\bbx_i.
\end{equation}
Upon defining the aggregate Lagrangian maximizer vector $\bbx(\bblambda):=[\bbx_1(\bblambda);\dots;\bbx_n(\bblambda)]$ as the concatenation of the local maximizers $\bbx_i(\bblambda)$, we can define the dual function as $h(\bblambda) := \ccalL(\bbx(\bblambda),\bblambda)$ and the dual problem as the minimization of the dual function,
\begin{align}\label{eq_dual_problem}
\bblambda^* &:= \underset{\bblambda}{\text{ argmin}} \ h(\bblambda) \\
& :=
\underset{\bblambda}{\text{ argmin}}
\sum_{i=1}^n f_i\left(\bbx_i(\bblambda)\right) +\!\! \sum_{(i,j) \in \ccalE} \!\!\bblambda_{ij}^T\left(\bbx_i(\bblambda) - \bbx_j(\bblambda)\right).\nonumber
\end{align}
For concave problems that satisfy minimal constraint qualifications that we assume to hold, the dual problem in \eqref{eq_dual_problem} is equivalent to the primal problem in \eqref{eq_primal_problem}. In particular, the optimal primal variable at node $i$ can be recovered as $\bbx_i(\bblambda^*)=\bbx^*$ if the optimal multiplier $\bblambda^*$ is known [cf.\eqref{eq_lagrangian_maximizers}].

An important feature of the dual function $h(\bblambda)$ is that its gradients can be computed locally as well. Specifically, it follows from the definition of the dual objective function $h(\bblambda)$ in \eqref{eq_dual_problem} that the partial derivative of $h(\bblambda)$ with respect to the component $\bblambda_{ij}$ can be written as the constraint slack
\begin{equation}\label{eq_dual_derivative}
    \bbg_{ij}(\bblambda) := \frac{\partial h(\bblambda)}{\partial \bblambda_{ij}} 
                    = \bbx_i(\bblambda) - \bbx_j  (\bblambda).
\end{equation}
Observe that to solve the program in \eqref{eq_lagrangian_maximizers}, node $i$ requires access to the local multipliers $\bblambda_{ij}$ and the dual variables $\bblambda_{ji}$ of its neighbors $j\in n_i$. Likewise, to evaluate the gradients $\bbg_{ij}(\bblambda)$ in \eqref{eq_dual_derivative}, node $i$ needs access to the local Lagrangian maximizer $\bbx_i  (\bblambda)$ and the neighboring maximizer $\bbx_j(\bblambda)$. It follows that gradient descent in the dual function can be implemented distributedly by relying on local operations and communication with neighboring nodes \cite{cRabbatNowak04}.

For future reference we emphasize that although the primal objective function $f$ is strongly concave, the dual objective function $h$ is not necessarily strongly convex. This fact requires the use of regularizations in the quasi-Newton algorithm that we will propose in Section \ref{sec_dbfgs}. We study this regularization in the following section.

%
\subsection{Regularized BFGS} \label{sec_bfgs_cent}

Dual gradient descent can be implemented in a decentralized manner and proven to converge to optimal arguments. However, convergence is slow when the condition number of the dual function is large. In this paper we propose a {\it decentralized} quasi-Newton method to overcome this limitation. In {\it centralized} settings, the idea of quasi-Newton methods is to alter the descent direction by premultiplying the dual gradient with an approximation of its Hessian inverse. Specifically, consider a time index $t$ and let $\bblambda(t)$ denote the dual variable iterate at time $t$ and $\bbg(t)=\bbg(\bblambda(t))$ be the corresponding gradient. Further, introduce a step size $\eps(t)$ and a symmetric positive definite matrix $\bbB(t)\in\reals^{m p \times m p}$ to define the dual quasi-Newton method through the recursion
\begin{equation}\label{eq_qnewt_update}
   \bblambda(t+1)  = \bblambda(t) - \epsilon(t)\bbB(t)^{-1}\bbg(t) 
              := \bblambda(t)+ \epsilon(t)\bbd(t) ,
\end{equation}
where we have defined the descent direction $\bbd(t):=-\bbB(t)^{-1}\bbg(t)$ in the second equality. If we substitute the matrix $\bbB(t)$ by the dual Hessian $\nabla^2h(\bblambda(t))$, we recover the update of Newton's method. The idea of quasi-Newton methods is that to design the matrix $\bbB(t)$ as an approximation of the the dual Hessian $\nabla^2h(\bblambda(t))$, while avoiding the cost of its evaluation. Various quasi-Newton methods are known to accomplish this feat, with the most common being the method of Broyden-Fletcher-Goldfarb-Shanno (BFGS) \cite{nocedal2006numerical}. To describe BFGS begin by defining the variable variation $\bbv(t)$ and the gradient variation $\bbr(t)$ as
\begin{equation}\label{eq_bfgs_vars} 
    \bbv(t) := \bblambda(t+1) - \bblambda(t), \quad 
    \bbr(t) := \bbg(t+1) - \bbg(t).
\end{equation}
The idea of {\it regularized} BFGS \cite{mokhtari2014res} is to find a matrix at each iteration $t+1$ that: (i) Satisfies the secant condition $\bbB(t+1)\bbv(t) =   \bbr(t)$. (ii) Is closest to the previous Hessian approximation matrix $\bbB(t)$ with respect to a differential entropy measure. (iii) Has a smallest eigenvalue not smaller than a pre-specified constant $\gamma$. To express this matrix introduce the modified gradient variation vector $\tbr(t)$ with regularization constant $\gamma>0$,
\begin{equation}\label{mod_grad_var}
    \tbr(t):=\bbg(t+1)-\bbg(t)-\gamma\bbv(t).
\end{equation}
The Hessian approximation $\bbB(t)$ of regularized BFGS can then be computed by recursive application of 
\begin{equation}\label{eq_reg_bfgs_cent} 
   \bbB(t+1) =   \bbB(t) 
               + \frac{\tbr(t) \tbr(t)^T}
                      {\tbr(t)^T \bbv(t)} 
               - \frac{\bbB(t) \bbv(t) \bbv(t)^T \bbB(t)}
                      {\bbv(t)^T \bbB(t) \bbv(t)}+\gamma\bbI.
\end{equation}
The matrix $\bbB(t+1)$ in \eqref{eq_reg_bfgs_cent} is the closest to $\bbB(t)$ in terms of relative entropy among all the matrices that satisfy the original secant condition and whose smallest eigenvalue is at least $\gamma$ \cite[Proposition 1]{mokhtari2014res}. We emphasize that the  use of the modified gradient variation in lieu of the regular gradient variation is necessary to maintain validity of the secant $\bbB(t+1)\bbv(t)=\bbr(t)$ condition which is the property that endows BFGS with appealing convergence traits; see \cite{mokhtari2014res} for details.

Both, the variable iteration in \eqref{eq_qnewt_update} and the matrix update in \eqref{eq_reg_bfgs_cent}, require centralized operations. In particular, to evaluate the inner product $\tbr(t)^T \bbv(t)$ in \eqref{eq_reg_bfgs_cent} or to compute the product $\bbB(t)^{-1}\bbg(t)$ in \eqref{eq_qnewt_update} nodes need access to global information. The goal of this paper is to introduce a variation of the regularized BFGS method that maintains the secant condition and is implementable in a decentralized manner.  

%
\section{Decentralized BFGS}\label{sec_dbfgs}

We propose a decentralized implementation of BFGS (D-BFGS), in which each node approximates the curvature of its local cost function and its neighbors. In doing so, each node computes and stores a local Hessian inverse approximation. As in regularized BFGS, an important feature of D-BFGS is that while it can be formulated in a decentralized manner, the global descent of the algorithm still satisfies the original secant condition. To study the details, recall $n_i$ as the neighborhood of node $i$ and define $m_i:=| n_i|$ as the number of neighbors of node $i$. We introduce the local dual vector $\bblambda_{i}(t)\in\reals^{m_i p}$ of node $i$ as the concatenation of the dual variables $\bblambda_{ij}(t)$ where $j\in n_i$. We use this definition to introduce the local variable variation at node $i$ as
\begin{equation}\label{local_variable_var}
\bbv_i(t):=\frac{\bblambda_i(t+1)-\bblambda_i(t)}{m_i+1}. 
\end{equation}
Nodes can exchange their local dual variable $\bblambda_i(t)$ with each other and construct a concatenated \textit{neighborhood} dual variable  $\bblambda_{ n_i}\in\reals^{M_ip}$ where $M_i=m_i + \sum_{j\in n_i} m_j$. The local variable variation in \eqref{local_variable_var} can subsequently be extended to  the neighborhood variable variation vector $\tbv_{ n_i}(t)\in\reals^{M_ip}$. Specifically, denote by $\bbD_{ n_i}\in\reals^{M_ip\times M_ip}$ the diagonal matrix such that its components corresponding to node $j$ are $1/(m_j+1)$. The modified neighborhood variable variation $\tbv_{ n_i}(t)$ of node $i$ is then given by
\begin{equation}\label{eq_dbfgs_vars} 
\tbv_{ n_i}(t) = \bbD_{ n_i} \left[ \bblambda_{ n_i}(t+1) - \bblambda_{ n_i}(t) \right].
\end{equation}
Likewise, we can define the local gradient $\bbg_i(t)\in \reals^{m_ip}$ at node $i$ as the concatenation of the partial derivates ${\partial h(\bblambda)}/{\partial \bblambda_{ij}} =\bbx^*_j(\bblambda) - \bbx^*_i  (\bblambda)$ for all $j\in n_i$. These can be then exchanged between neighboring nodes to construct a neighborhood gradient $\bbg_{ n_i}(t)\in\reals^{M_ip}$ and modified neighborhood gradient variation $\tbr_{ n_i}(t)\in\reals^{M_ip}$. With a  regularization constant $\gamma > 0$, the modified neighborhood gradient variation $\tbr_{ n_i}(t)$ is given by
\begin{equation}\label{eq_dbfgs_grads} 
\tbr_{ n_i}(t) =\bbg_{ n_i}(t+1) - \bbg_{ n_i}(t)  -\gamma \tbv_{ n_i}(t).
\end{equation}
%

Thus, we have defined the neighborhood variable variation $\tbv_{ n_i}(t)$ and modified gradient variation $\tbr_{ n_i}(t)$ so that they are suitable for decentralized settings. We introduce $\bbB^i(t)$ as the neighborhood Hessian approximation matrix, which is updated as 
\begin{align}\label{dbfgs}
\bbB^i(t+1)&= \bbB^i(t) + \frac{\tbr_{ n_i}(t) \tbr_{ n_i}(t)^T}{\tbr_{ n_i}(t)^T \tbv_{ n_i}(t)}   \\
&\quad- \frac{\bbB^i(t) \tbv_{ n_i}(t) \tbv_{ n_i}(t)^T \bbB^i(t)}{\tbv_{ n_i}(t)^T \bbB^i(t) \tbv_{ n_i}(t)} + \gamma \bbI.\nonumber
\end{align}
The update in \eqref{dbfgs} differs from \eqref{eq_reg_bfgs_cent} in its use of neighborhood variable and modified gradient variation as defined in \eqref{eq_dbfgs_vars} and \eqref{eq_dbfgs_grads}, respectively, instead of the variation vectors in \eqref{eq_bfgs_vars} and \eqref{mod_grad_var}. 
%
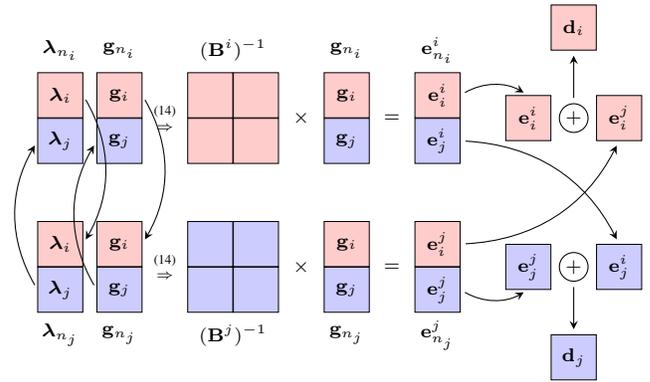
\begin{figure}\centering

\def \thisplotscale {0.6}
\def \unit {\thisplotscale cm}

\tikzstyle {block}        = [draw, very thin,
                             rectangle, 
                             minimum height = 1*\unit,
                             minimum width  = 1*\unit]                       

\tikzstyle {blue block}   = [block,
                             fill = blue!20]

\tikzstyle {green  block} = [block,
                             fill = green!20]

\tikzstyle {red block}    = [block,
                             fill = red!20]

\def \deltalabel{ 0.5}

{\fontsize{7}{7}\selectfont\begin{tikzpicture}[x = 1*\unit, y=1*\unit, 
                          shorten >=2pt, shorten <=2pt]

%
%
\path (0,0)        node [red block       ] (Bine) {};
\path (Bine.south) node [red block, below] (Bise) {};
\path (Bise.west)  node [red block, left ] (Bisw) {};
\path (Bisw.north) node [red block, above] (Binw) {};
\path (Bine.north west) ++ (0, \deltalabel) node {$(\bbB^i)^{-1}$}; 
%
%
\path (Binw.west) ++ (-1,0)   node [red   block, left] (rnin) {$\bbg_i$};
\path (Bisw.west) ++ (-1,0)   node [blue  block, left] (rnis) {$\bbg_j$};
\path (rnin.west) ++ (-0.3,0) node [red   block, left] (vnin) {$\bblambda_i$};
\path (rnis.west) ++ (-0.3,0) node [blue  block, left] (vnis) {$\bblambda_j$};
\path (rnin.north) ++ (0, \deltalabel) node {$\bbg_{n_i}$}; 
\path (vnin.north) ++ (0, \deltalabel) node {$\bblambda_{n_i}$}; 
%
%
\path (Bine.east) ++ (1,0) node [red   block, right] (gnin) {$\bbg_i$};
\path (Bise.east) ++ (1,0) node [blue  block, right] (gnis) {$\bbg_j$};
\path (gnin.north) ++ (0, \deltalabel) node {$\bbg_{n_i}$}; 
%
%
\path (gnin.east) ++ (1,0) node [red   block, right] (dnin) {$\bbe_i^i$};
\path (gnis.east) ++ (1,0) node [blue  block, right] (dnis) {$\bbe_j^i$};
\path (dnin.north) ++ (0, \deltalabel) node {$\bbe^i_{n_i}$}; 
%
%
\path (dnin.south east) ++ (1,0) node [red block, right] (dil) {$\bbe_i^i$};
\path (dil.east)        ++ (1,0) node [red block, right] (dir) {$\bbe_i^j$};
\path (rnin.south) -- (Binw.south)  node [midway] {$\text{\eqref{dbfgs}} \atop
                                                    \Rightarrow$};
\path (Bine.south) -- (gnin.south)  node [midway] {$\times$};
\path (gnin.south) -- (dnin.south)  node [midway] {$=$};
\path (dil.east)   -- (dir.west)    node [midway, draw, circle, inner sep=1] (plus) {$+$};
\path[draw, -stealth] (plus) -- ++(0,1.5) node [red block, above] {$\bbd_i$};;

%
%
%
\path (0, -3.3)    node [blue block       ] (Bjne) {};
\path (Bjne.south) node [blue block, below] (Bjse) {};
\path (Bjse.west)  node [blue block, left ] (Bjsw) {};
\path (Bjsw.north) node [blue block, above] (Bjnw) {};
\path (Bjse.south west) ++ (0,-\deltalabel) node {$(\bbB^j)^{-1}$}; 
%
%
\path (Bjnw.west) ++ (-1,0)   node [red   block, left] (rnjn) {$\bbg_i$};
\path (Bjsw.west) ++ (-1,0)   node [blue  block, left] (rnjs) {$\bbg_j$};
\path (rnjn.west) ++ (-0.3,0) node [red   block, left] (vnjn) {$\bblambda_i$};
\path (rnjs.west) ++ (-0.3,0) node [blue  block, left] (vnjs) {$\bblambda_j$};
\path (rnjs.south) ++ (0,-\deltalabel) node {$\bbg_{n_j}$}; 
\path (vnjs.south) ++ (0,-\deltalabel) node {$\bblambda_{n_j}$}; 
%
%
\path (Bjne.east) ++ (1,0) node [red   block, right] (gnjn) {$\bbg_i$};
\path (Bjse.east) ++ (1,0) node [blue  block, right] (gnjs) {$\bbg_j$};
\path (gnjs.south) ++ (0,-\deltalabel) node {$\bbg_{n_j}$}; 
%
%
\path (gnjn.east) ++ (1,0) node [red   block, right] (dnjn) {$\bbe_i^j$};
\path (gnjs.east) ++ (1,0) node [blue  block, right] (dnjs) {$\bbe_j^j$};
\path (dnjs.south) ++ (0,-\deltalabel) node {$\bbe^j_{n_j}$}; 
%
%
\path (dnjn.south east) ++ (1,0) node [blue block, right] (djl) {$\bbe_j^j$};
\path (djl.east)        ++ (1,0) node [blue block, right] (djr) {$\bbe_j^i$};
%
%
\path (rnjn.south) -- (Bjnw.south)  node [midway] {$\text{\eqref{dbfgs}} \atop
                                                    \Rightarrow$};
\path (Bjne.south) -- (gnjn.south)  node [midway] {$\times$};
\path (gnjn.south) -- (dnjn.south)  node [midway] {$=$};
\path (djl.east)   -- (djr.west)    node [midway, draw, circle, inner sep=1] (plus) {$+$};
\path[draw, -stealth] (plus) -- ++(0,-1.5) node [blue block, below] {$\bbd_j$};;

%
%
\path (dnin.east) edge [-stealth, black!99, bend left]  (dil.north);
\path (dnis.east) edge [-stealth, black!99, bend left]  (djr.north);
\path (dnjs.east) edge [-stealth, black!99, bend right] (djl.south);
\path (dnjn.east) edge [-stealth, black!99, bend right] (dir.south);
%
%
%
%
\path (rnin.east) edge [-stealth, black!99, bend left] (rnjn.east);
\path (rnjs.west) edge [-stealth, black!99, bend left] (rnis.west);
\path (vnin.east) edge [-stealth, black!99, bend left] (vnjn.east);
\path (vnjs.west) edge [-stealth, black!99, bend left] (vnis.west);

\end{tikzpicture}}
\caption{D-BFGS variable flow. Nodes exchange variable and gradient variations -- $\bblambda_i$ and $\bbg_i$ sent to $j$ and $\bblambda_j$ and $\bbg_j$ sent to $i$ -- that they use to determine local curvature matrices -- $\bbB^i$ and $\bbB^j$. They then use exchanged gradients to compute descent directions -- $\bbe_{n_i}^i$ and  $\bbe_{n_j}^j$. These contain a piece to add locally -- $\bbe_i^i$ stays at node $i$ and $\bbe_j^j$ stays at node -- and a piece to add at neighbors -- $\bbe_j^i$ is sent to node $j$ and $\bbe_i^j$ is sent to node $i$.}
\label{fig_variable_flow_diagram} \end{figure}

We define the descent direction of D-BFGS with normalization constant $\Gamma >0$ evaluated at node $i$ as
\begin{equation}\label{eq_direction_local}
\bbe^i_{ n_i}(t) := - (\bbB^i(t)^{-1} +\Gamma\bbD_{n_i})\ \bbg_{ n_i}(t).\end{equation}
Note that  $\Gamma\bbD_{n_i}$ is added to the Hessian inverse approximation $\bbB^i(t)^{-1}$ to ensure that the descent direction $\bbe^i_{ n_i}(t) $ is not null. 

The variable flow is demonstrated in Figure \ref{fig_variable_flow_diagram}.  Nodes exchange variable and gradient information to compute local Hessian inverse approximations $\bbB^i(t)^{-1}$ and neighborhood descent direction $\bbe^i_{n_i}(t)$. The descent direction $\bbe^i_{ n_i}(t)$ contains a descent direction for node $i$ and its neighbors $j \in  n_i$. Nodes therefore exchange with their neighbors the parts of their locally computed descent direction pertaining to them. To be more precise, denote $\bbe_{ j}^i(t) \in \mathbb{R}^{m_jp} = [\bbe^i_{ n_i}(t)]_j$ as the component of the descent direction $\bbe^i_{ n_i}(t)$ evaluated at node $i$ that belongs to the neighbor $j$. Node $i$ computes its full descent direction $\bbd_i(t)$ as the sum of locally computed descent directions $\bbe_{i}^i(t)$ and the parts received from its neighbors $\bbe^j_{i}(t)$, i.e., 
\begin{equation}\label{eq_direction_coord}
\bbd_{i}(t) := \bbe^i_{i}(t)+ \sum_{j \in  n_i}  \bbe^j_{i}(t).
\end{equation}
The local variable $\bblambda_i$ at node $i$ is then updated using the full descent direction $\bbd_{i}(t)$ by
\begin{equation}
\bblambda_{i}(t+1) = \bblambda_{i}(t) +\epsilon(t) \bbd_i(t).
\label{eq_descent_coord}
\end{equation}
Note that at each step $t$ the primal variable $\bbx_i(t+1)$ can subsequently be recovered as the Lagrangian maximizer with respect to $\bblambda(t+1)$, i.e.
\begin{equation}\label{eq_primal_update}
\bbx_i(t+1)=\bbx_i(\bblambda(t+1)).
\end{equation}

%
\begin{algorithm}[t] 
\setstretch{1.35}
\small{\begin{algorithmic}[1]

  \REQUIRE
  $\bbB^i(0), \bblambda_{i}(0),\bbg_i(0), \bblambda_{n_i}(0), \bbg_{n_i}(0)$
  \FOR{$t = 0,1,2, \hdots$}    
  \STATE Compute $\bbe^i_{ n_i}(t)\! =\! - (\bbB^i(t)^{-1} \!+\!\Gamma\bbD_{n_i}) \bbg_{ n_i}(t)$ [cf. \eqref{eq_direction_local}]
  \STATE Exchange $\bbe^i_{j}(t)$ with neighbors $j \in n_i$
  \STATE Compute descent dir. $\displaystyle{\bbd_{i}(t) := \bbe^i_{i}(t)+ \sum_{j \in  n_i}  \bbe^j_{i}(t).}$
  \STATE Update local variable $\bblambda_{i}(t+1) = \bblambda_{i}(t) +\epsilon(t) \bbd_i(t)$ [cf.\eqref{eq_descent_coord}] and exchange with neighbors
      \STATE Compute $\bbx_i(\bblambda(t+1))$ and exchange with neighbors  [cf. \eqref{eq_lagrangian_maximizers}]\\
       \STATE Compute $\bbg_{ij}(t+1) $ and exchange with neighbors
 [cf. \eqref{eq_dual_derivative}]\\ $\displaystyle{\bbg_{ij}(t\!+\!1)= \bbx_i(\bblambda(t\!+\!1)) - \bbx_j  (\bblambda(t+1))}$ 
  \STATE Compute $\tbv_{ n_i}(t),\tbr_{ n_i}(t),\bbB^i(t+1)$ [cf.\eqref{eq_dbfgs_vars}--\eqref{dbfgs}] 
  \ENDFOR
\end{algorithmic}}
\caption{D-BFGS method at node $i$}
\label{alg_dbfgs}
\end{algorithm}

The summary of the algorithm performed for node $i$ is outlined in Algorithm \ref{alg_dbfgs}. Each node begins with an initial dual variable $\bbv_i(0)$, and Hessian approximation $\bbB^i(0)$, and gradient $\bbg_i(0)$. Nodes initially exchange local variable and gradients to construct initial neighborhood variable $\bblambda_{n_i}(0)$ and gradient $\bbg_{n_i}(0)$. For each step $t$, nodes compute their neighborhood descent direction $\bbe^i_{ n_i}(t) $ in Step 3 and exchange the  descent elements $\bbe^j_{ n_i}(t) $ with their neighbors in Step 3 to compute the full descent direction $\bbd_i(t) $ as in Step 4. They use the full descent direction $\bbd_i(t) $ to update the variable $\bblambda_i(t+1)$ and exchange it with their neighbors to form $\bblambda_{n_i}(t+1)$ in Step 5. Then, they use their neighbor variables $\bblambda_j(t+1)$ to compute an updated Lagrangian maximizer $\bbx_{i}(\bblambda(t+1))$ in Step 6, which it then exchanges with its neighbors.
With access to the Lagrangian maximizer $\bbx_{i}(\bblambda(t+1))$ of their neighbors, the gradient $\bbg_{ij}(t+1)$ can be updated and then exchanged as in Step 7. In Step 8, nodes compute their modified variable $\tbv_{n_i}(t)$ and gradient $\tbr_{n_i}(t)$ variations that are required for computing the updated neighborhood Hessian approximation matrix $\bbB^i(t+1)$.

\begin{remark}\normalfont
We stress here the need in using the normalization matrix $\bbD_{n_i}$ used in \eqref{eq_dbfgs_vars} for defining the variable variation $\bbv_{n_i}(t)$. Including $\bbD_{n_i}$ in the definition ensures that the global Hessian inverse approximation matrix satisfies the global secant condition. To be more precise, consider the regularized neighborhood Hessian inverse approximation $\bbB^i(t)^{-1} +\Gamma\bbI\in\reals^{M_ip\times M_ip}$. We define $\bbH^i(t) \in \mathbb{R}^{mp \times mp}$ to be the block sparse matrix with respect to $n_i$ that has a dense sub-matrix $\bbB^i(t)^{-1}$. Further define $\hbD_{n_i} \in \reals^{mp \times mp}$ to be the black sparse matrix with respect to $n_i$ that has a dense submatrix of $\bbD_{n_i}$. Considering this definition the descent direction $\bbd(t)$ for the global concatenated variable vector $\bblambda(t)$ can be written as 
\begin{align}\label{eq_direction_dist}
\bbd(t) &:=  -\sum_{i=1}^n \left[ \bbH^i(t) + \Gamma \hbD_{n_i} \right]  \bbg(t) \nonumber \\
&:= - \left[ \bbH(t) + \Gamma \bbI \right] \bbg(t).
\end{align}
The expression in \eqref{eq_direction_dist} states that the matrix $\bbH(t) : = \sum_{i=1}^n \bbH^i(t)$ is the global Hessian inverse approximation. It is easily verifiable that the matrix $\bbH(t)$ satisfies the global secant condition, i.e., $\bbv(t) = \bbH(t) \bbr(t)$ which justifies the normalization of the variable variation in \eqref{eq_dbfgs_vars}.
\end{remark}

\begin{remark}\label{remark_inner_product}\normalfont
As in the case of centralized BFGS, it is necessary that the neighborhood inner product $\tbr_{ n_i}(t)^T \tbv_{ n_i}(t)$ be positive in order for $\bbB^i(t)$ to be well defined. In the decentralized case, however, due to the truncating of the gradient and variable vectors we cannot guarantee that  this inner product is positive even when the functions are strongly convex. As such, in practice the condition $\tbr_{ n_i}(t)^T \tbv_{ n_i}(t)>0$ must be verified at every iteration, otherwise we do not update the Hessian approximation matrix, i.e. $ \bbB^i(t+1) = \bbB^i(t)$. Note that for these iterations the local secant condition for node $i$ is not satisfied.
\end{remark}

%
\subsection{Asynchronous implementation} \label{sec_async_dbfgs}
Given the coordination and communication cost required to implement D-BFGS in Algorithm \ref{alg_dbfgs}, we also consider the D-BFGS algorithm in the asynchronous setting. Our model for asynchronicity follows that used in \cite{bertsekas1989parallel}, in which nodes perform computations and communications out of sync with their neighbors. Consider that the time indices are partitioned so that node $i$'s primary computation, namely the computation of descent direction $\bbe^i_{n_i}(t)$, requires multiple time iterates to complete. For each node $i$, we define a set $T^i \subseteq \mathbb{Z}^+$ of all time indices in which node $i$ is available to send and receive information.

We further define two functions that specify the asynchronicity between nodes. For each node $i$, we define a function $\pi^i(t)$ that returns the most recent time node prior to $t$ node $i$ was available, i.e.
\begin{equation}
\pi^i(t) := \text{max} \{\hat{t} \mid \hat{t} < t, \hat{t} \in T^i\} \label{eq_time_function}.
\end{equation}
Moreover, we define a function $\pi^i_j(t)$ that returns the most recent time node $j$ sent information that has been received by node $i$ by time $t$, or explicitly, 
\begin{equation}
 \pi_j^i(t) := \pi^j(\pi^i(t)) \label{eq_time_function_neighbor}.
 \end{equation}
In the asynchronous setting, the superscript notation used to denote locally computed information now additionally signifies a node's current knowledge its neighbors, i.e.
 \begin{align}
 \bblambda_{j}^i(t) &:= \bblambda_j(\pi_j^i(t)) \\
\bblambda_{n_i}^i(t) &=[ \bblambda_j^i(t)]_{j \in n_i}.
 \label{eq_local_var_async}
 \end{align}
 We emphasize that $\bblambda_{j}^i(t) \neq \bblambda_j^k(t)$ for any two nodes $i$ and $k$ at any time $t$. We consider as the current global variable state $\bblambda(t)$ the concatenation of each node's current knowledge of its own variable, i.e. $\bblambda(t) := [\bblambda^i_i(t); \hdots; \bblambda^n_n(t)]$. We subsequently use the same notation for local gradients $\bbg^i_j(t)$ and descent directions $\bbe^i_j(t)$.

 We assume at any time $t \in T^i$ that node $i$ does three things: (i) It reads the variable, gradient, and descent directions from neighboring nodes $j \in n_i$ sent while it was busy. (ii) It updates its local variables and gradient using the descent direction is has just finished computing as well as the descent directions it has received from its neighbors. (iii) Node $i$ can send its locally computed descent direction as well as its updated variable and gradient info. To state in more explicit terms, node $i$ performs the following update to its own variable at all times,
\begin{align}\label{eq_update_local_async}
    \bblambda^i_{i}(t+1) = \bblambda^i_{i}(t) + \epsilon(t) \bbd_{i}(t),
\end{align} 
where $\bbd_{i}(t)$ is the decent for the $i$th variable $\bbx_{i}(t)$ at $t$,
\begin{align}
 \bbd_{i}(t) &=  \begin{cases}  \bbe^i_{i}(t) +  \sum_{j \in n_i}  \bbe^j_{i}(t)  & \text{if $t \in T^i$} \\
       \bb0. & \text{otherwise}.
       \end{cases}
       \label{eq_descent_local}
 \end{align}
If $t \in T^i$, node $i$ applies all descent directions available, otherwise it does nothing. Observe that the descent direction in \eqref{eq_descent_local} contains descents calculated with information from time $\pi^i(t)$ as well as the times $\pi^i_j(t)$ that neighbor $j$ most recently updated its local variable. 
 
To specify the asynchronous version of the D-BFGS algorithm, we reformulate the variable and gradient differences, $\tbv_{n_i}^i(t)$ and $\tbr_{n_i}^i(t)$ for the asynchronous case:
\begin{align}
\tbv_{n_i}^i(t) &= \bbD_{n_i} \left[ \bblambda_{n_i}^i(t+1) - \bblambda_{n_i}^i(t) \right], \label{eq_dbfgs_vars_async} \\
\tbr_{n_i}^i(t) &= \bbg_{n_i}^i(t+1) - \bbg_{n_i}^i(t) - \gamma \tbv_{n_i}^i(t). \label{eq_dbfgs_grads_async}
\end{align}
The computation of the local asynchronous update matrix $\bbB^i(t)$ and the corresponding descent direction $\bbe^i_{n_i}(t)$ follows respectively \eqref{dbfgs} and \eqref{eq_direction_local} exactly as in the synchronous setting, now using asynchronous $\tbv_{n_i}^i(t)$ and $\tbr_{n_i}^i(t)$ in place of $\tbv_{n_i}(t)$ and $\tbr_{n_i}(t)$, respectively.

The complete asynchronous algorithm is outlined in Algorithm \ref{alg_async_dbfgs}. Each node begins with an initial variable $\bblambda_i(0)$, Hessian approximation $\bbB^i(0)$, gradient $\bbg_i(0)$, and descent component $\bbe^i_i(0)$. At each time index $t$, they begin by reading the variables of neighbors $\bbe^j_{i}(t), \bblambda^i_{j}(t), \bbg^j_{i}(t)$ in Step 2 and construct neighborhood variables. The aggregated descent direction $\bbd_i(t) $ is used to update variables $\bblambda_i(t+1)$, $\bbx(\bblambda(t+1))$, and $\bbg_i(t+1)$ in Step 3. Then, with the updated local variable $\bblambda_i(t+1)$ and gradient $\bbg_i(t+1)$, the node computes the D-BFGS variables $\tbv^i_{n_i}(t)$, $\tbr^i_{n_i}(t)$, and $\bbB^i(t+1)$ in Step 4. In Step 5, node $i$ computes the next descent direction $\bbd^i_{n_i}(t+1)$, and sends its variables to neighbors in Step 6.

%
\begin{algorithm}[t]
\setstretch{1.35}
{\small\begin{algorithmic}[1]
  \REQUIRE $\bbB^i(0)$, 
           $\bblambda_{i}(0)$, 
           $\bbg_{i}(0)$, 
           $\bbd^i_{n_i}(0)$ [cf. \eqref{eq_direction_local}]
  \FOR{$t \in T^i$}
  \STATE Read $\bbe^j_{i}(t), \bblambda^i_{j}(t), \bbg^j_{i}(t)$ from neighbors $j \in n_i$
  \STATE Update $\bblambda_{i}(t+1), \bbx(\bblambda(t)), \bbg_{i}(t+1)$ [cf. \eqref{eq_update_local_async}, \eqref{eq_descent_local}]
  \STATE Compute $\tbv^i_{n_i}(t),\tbr^i_{n_i}(t),\bbB^i(t+1)$ [cf. \eqref{eq_dbfgs_vars_async}, \eqref{eq_dbfgs_grads_async}, \eqref{dbfgs}]
  \STATE Compute $\bbe^i_{n_i}(t+1)$ [cf. \eqref{eq_direction_local}]
  \STATE Send $\bblambda_{i}(t+1), \bbg_{i}(t+1)$, $\bbe^i_{j}(t+1)$ to neighbors $j \in n_i$  

  \ENDFOR
\end{algorithmic}}
\caption{Asynchronous D-BFGS method at node $i$}
\label{alg_async_dbfgs}
\end{algorithm}

While Algorithm \ref{alg_async_dbfgs} follows a similar structure to the synchronous Algorithm \ref{alg_dbfgs}, we highlight that n the synchronous algorithm, coordination of four rounds of communication were required at each iteration of Algorithm \ref{alg_dbfgs} to properly communicate the dual variable, primal variable, and dual gradient information. In the asynchronous setting, only a single round of communication is possible at each time iteration, thus the communication burden is indeed reduced and the coordination between nodes rendered unnecessary. 

%
\section{Convergence Analysis}\label{sec_convergence}

In this section, we study the convergence properties of the D-BFGS method and show that the sequence of primal iterates $\bbx_i$ generated by D-BFGS converges to the optimal argument $\tbx^*$ of \eqref{original_problem}. In proving these results we make the following assumptions.

\begin{assumption} \label{assumption1}
The local objective functions $f_i(\bbx)$ are differentiable and strongly convex with parameter $\mu>0$.
\end{assumption} 

The strong convexity of local functions $f_i$ implies that the aggregate function $f(\bbx)=\sum_{i=1}^n f_i(\bbx_i)$ is also strongly convex with constant $\mu$. Define the oriented incidence matrix $\bbA\in\reals^{mp\times n}$ where the $A_{ev}$ component is $1$ if edge $e$ is started from node $v$ and is $-1$ if edge $e$ is ended at node $v$; otherwise $A_{ev}=0$. It is not difficult to see that the Hessian of the dual function can be written as $\nabla^2 h(\bblambda) = \bbA \big(\nabla^2 f(\bbx(\bblambda))\big)^{-1} \bbA^T$ from where we conclude that the eigenvalues of the dual function Hessian are upper bounded by $4n/\mu$ (note: $\bbA^T \bbA$ is two times the graph Laplacian matrix). In turn, this implies that the dual function gradients $\bbg(\bblambda)$ are Lipschitz continuous with constant $4n/\mu$,
 \begin{equation}\label{lip_dual_grad}
 \|\bbg(\bblambda) - \bbg(\tblambda)\| \leq  \frac{4n}{\mu}\|\bblambda - \tblambda\|.
 \end{equation}
We additionally make a further assumption regarding the inner product of neighborhood variable and gradient variations.
 \begin{assumption} \label{assumption2}
For all $i$ and $t$, the inner product between the neighborhood modified variable and gradient vector variations is strictly positive, i.e. $\tbv_{n_i}^T \tbr_{n_i} > 0$.
\end{assumption} 
This assumption is necessary to ensure all local Hessian approximations are well defined in \eqref{dbfgs}. While this assumption does not always hold in practice, we use it regardless to simplify analysis. We stress that, in the case the assumption is violated, setting $\bbB^i(t+1) = \bbB^i(t)$ (See Remark \ref{remark_inner_product}) does not have any bearing on the proceeding analysis.
  
We specify $\bbH(t) = \sum_{i=1}^n \bbH^i(t)$ as the  global Hessian inverse approximation and  $\bbd(t) = -[\bbH(t) + \Gamma \bbI] \bbg(t)$ as the global descent direction. The following lemma establishes the positive definiteness of the global descent direction.

\begin{lemma}\label{lemma_eigen_bound}
Consider the D-BFGS method introduced in \eqref{local_variable_var}-\eqref{eq_direction_dist}. Further, recall both the positive constants $\gamma$ and $\Gamma$ as the regularization parameters of D-BFGS and the definition of the global Hessian inverse approximation $\bbH(t) + \Gamma \bbI = \sum_{i=1}^n [ \bbH^i(t) + \Gamma \bbD_{n_i} ]$. The eigenvalues of the global Hessian inverse approximation $\bbH(t)$ are uniformly bounded as 
\begin{align}\label{eq_eigen_bounds_prop}
&\Gamma\bbI \preceq \bbH(t) + \Gamma \bbI \preceq \Delta \bbI :=  \left( \Gamma + \frac{n}{\gamma} \right) \bbI, 
\end{align} 
where $n$ is the size of network. 
\end{lemma}

 \begin{myproof}
The lower bound on $\bbH(t) + \Gamma \bbI$ follows immediately from the fact that $\bbH(t)$ is a sum of positive semidefinite matrices and is therefore a positive semidefinite matrix with eigenvalues greater than or equal to 0. The upper bound subsequently follows from the fact that each $\bbH_i(t)$ have eigenvalues upper bounded by $1/\gamma$, as the dense submatrix $\bbB^i(t)^{-1} \preceq 1/\gamma \bbI$. Then, the sum of $n$ such matrices recovers the upper bound in \eqref{eq_eigen_bounds_prop}.
 \end{myproof} $\newline$


The result in Lemma \ref{lemma_eigen_bound} shows that the eigenvalues of the global Hessian inverse approximation $\bbH(t)$ are uniformly bounded. Thus, D-BFGS descent direction $\bbd(t)\!=\!-\bbH(t)\bbg(t)$ is a valid descent direction. We use this result and convexity of the dual function $h$ to show that the sequence of dual objective function errors $h(\bblambda)-h(\bblambda^*)$ converges to null.

\begin{theorem}\label{theorem_convergence}
Consider the D-BFGS method introduced in \eqref{local_variable_var}-\eqref{eq_direction_dist}. If Assumption \ref{assumption1} holds true and $\epsilon(t)$ is chosen such that $\epsilon(t) < \Gamma \mu / (n \Delta^2)$, then the dual objective function error $h(\bblambda(t))-h(\bblambda^*)$ converges to zero at least in the order of $o(1/t)$, i.e.,
\begin{equation} \label{eq_convergence_result}
 h(\bblambda(t)) - h(\bblambda^*)\leq o\left(\frac{1}{t}\right).
\end{equation}
\end{theorem}


The proof for this result is standard, and is a small variation of that for gradient descent found in \cite[Proposition 1.3.3]{bertsekas1999nonlinear}. Theorem \ref{theorem_convergence} shows that the sequence of the dual objective function $h(\bblambda)$ generated by D-BFGS converges to the optimal dual function value $h(\bblambda^*)$. We conclude with a corollary establishing the convergence of the primal function error and the primal variables of the original problem in \eqref{eq_primal_problem}.
\begin{corollary} \label{cor_primal_convergence}
The sequence of primal function value error generated by the D-BFGS algorithm converges to zero, i.e.
\begin{equation}
\lim_{t\to \infty} f(\bbx^*) -  f(\bbx(t)) = 0.
\end{equation}
Furthermore, the sequence of primal variables $\bbx(t)$ converges to the optimal primal variable $\bbx^*$ at least in the order of $o(1/\sqrt{t})$, i.e.
\begin{equation} \label{eq_primal_convergence_result}
\| \bbx(t) - \bbx^* \| \leq o\left(\frac{1}{\sqrt{t}}\right)
\end{equation}
\end{corollary}

\begin{myproof}
See Appendix \ref{sec_cor_primal_convergence}.
\end{myproof}

We subsequently further provide a convergence result for the asynchronous implementation of DBFGS. The result uses the common assumption of partial asynchronicity, i.e. any two nodes are no more than some constant $B$ out of sync. We demonstrate convergence in the following theorem.

\begin{theorem}\label{theorem_convergence_asyn}
Consider the asynchronous D-BFGS algorithm proposed in \eqref{eq_update_local_async}-\eqref{eq_dbfgs_grads_async} and \eqref{dbfgs}-\eqref{eq_direction_coord}. If Assumptions \ref{assumption1} holds and the following partial asynchronicity property is satisfied for some $B > 0$, 
\begin{equation}
\max\{0, t-B+1\} \leq \pi^i_j(t) \leq t, \quad \text{for all }i,j,t,
\end{equation}
then there exists a stepsize $\epsilon(t)$ such that $\lim_{t \rightarrow \infty} \bbg(t) = 0$.
\end{theorem}
The proof for this result is omitted for space considerations and can be found in \cite[Theorem 3]{DBFGS2}. This theorem demonstrates that the uncoordinated and asynchronous implementation of DBFGS is still guaranteed to converge.

%
\section{Numerical Results} \label{sec_numerical_results}

We provide numerical results of the performance of D-BFGS and first order methods ADMM\cite{Schizas2008-1} and dual descent (DD) in solving a quadratic program. Consider the problem
\begin{align}\label{eq_simulation_problem}
\bbx^* := \argmax_{\bbx\in\reals^p}\sum_{i=1}^n -\frac{1}{2} \bbx^T \bbA_i \bbx - \bbb_i^T \bbx,
\end{align}
where $\bbA_i \in \mathbb{R}^{p \times p}$ is the positive definite matrix  and $\bbb_i \in \mathbb{R}^p$ is a random vector which are both available only at node $i$. In this case the Lagrangian maximizer and primal update in \eqref{eq_lagrangian_maximizers} can be computed with a closed form solution as
\begin{equation} \label{eq_primal_update_quadratic}
\bbx_i(\bblambda)= \bbA_i^{-1} \sum_{j\in  n_i}(\bblambda_{ji}-\bblambda_{ij}) +  \bbA_i^{-1}\bbb_i.
\end{equation}
The rest of the D-BFGS updates for this problem follow from \eqref{eq_dual_derivative} and \eqref{local_variable_var}-\eqref{eq_direction_dist}.

\begin{figure}[t]
\centering
\includegraphics[height=.25\textheight,width=\linewidth,keepaspectratio]{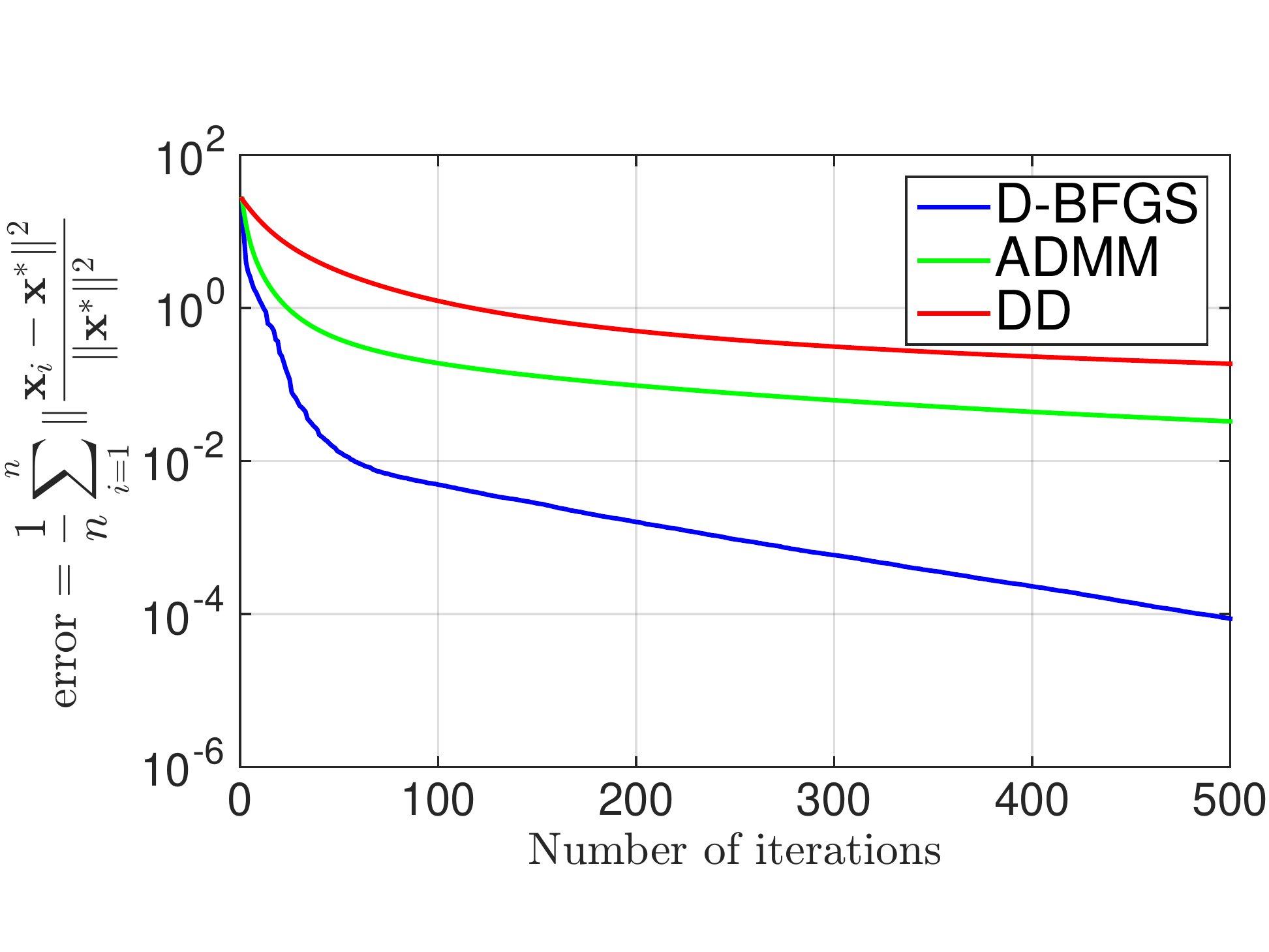}
\caption{Convergence path of the average distance to optimal primal variable vs. number of iterations for D-BFGS, ADMM, and DD for quadratic problem with condition number $10^2$.}\label{figure_simulation}
\end{figure}

\begin{figure*}	
	\centering
	\begin{subfigure}[t]{.3\textwidth}
		\centering
		\includegraphics[height=.13\textheight,width=\textwidth]{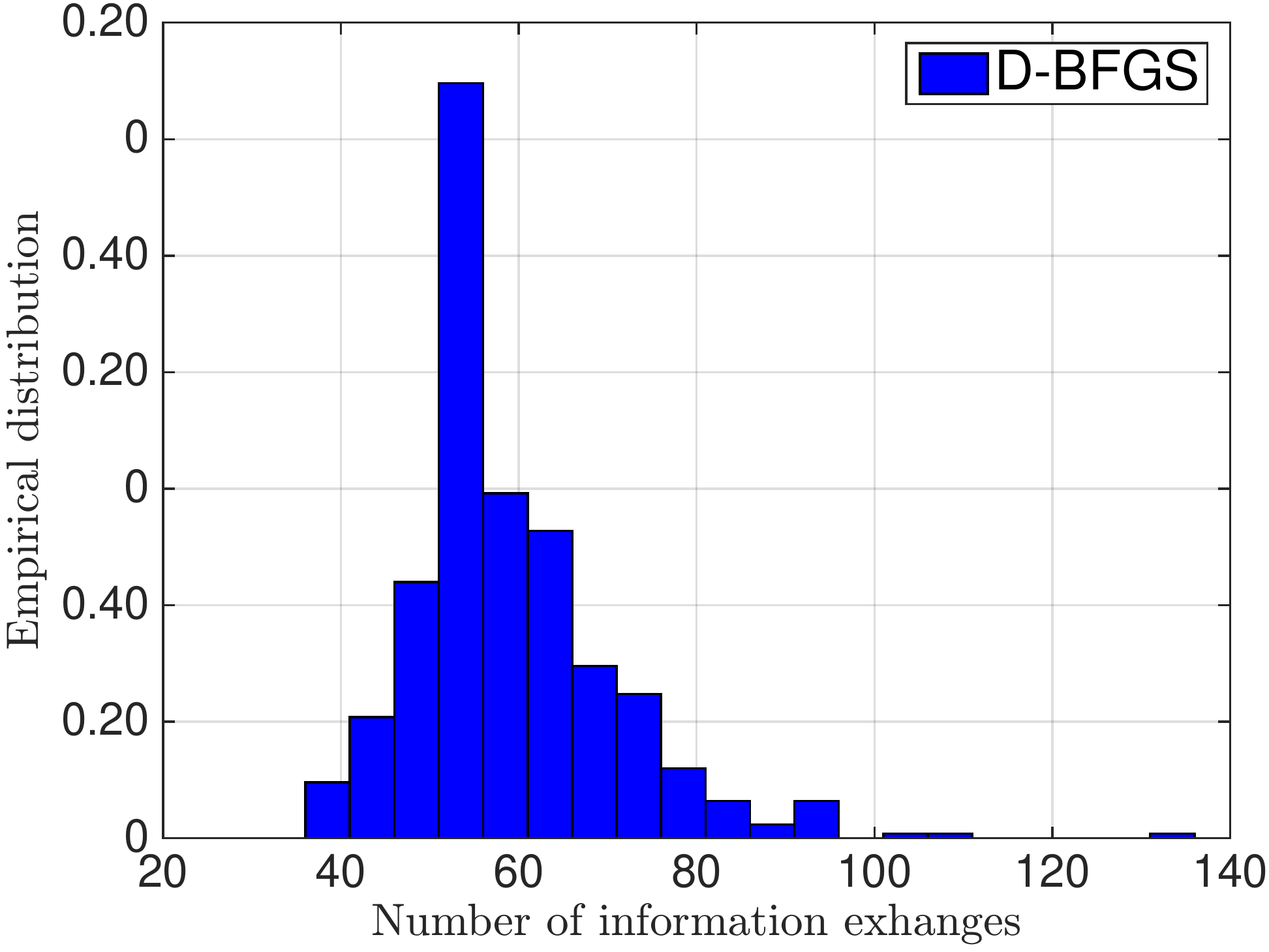}
		\caption{}\label{fig:1a}		
	\end{subfigure}
	\quad
	\begin{subfigure}[t]{.3\textwidth}
		\centering
		\includegraphics[height=.13\textheight,width=\textwidth]{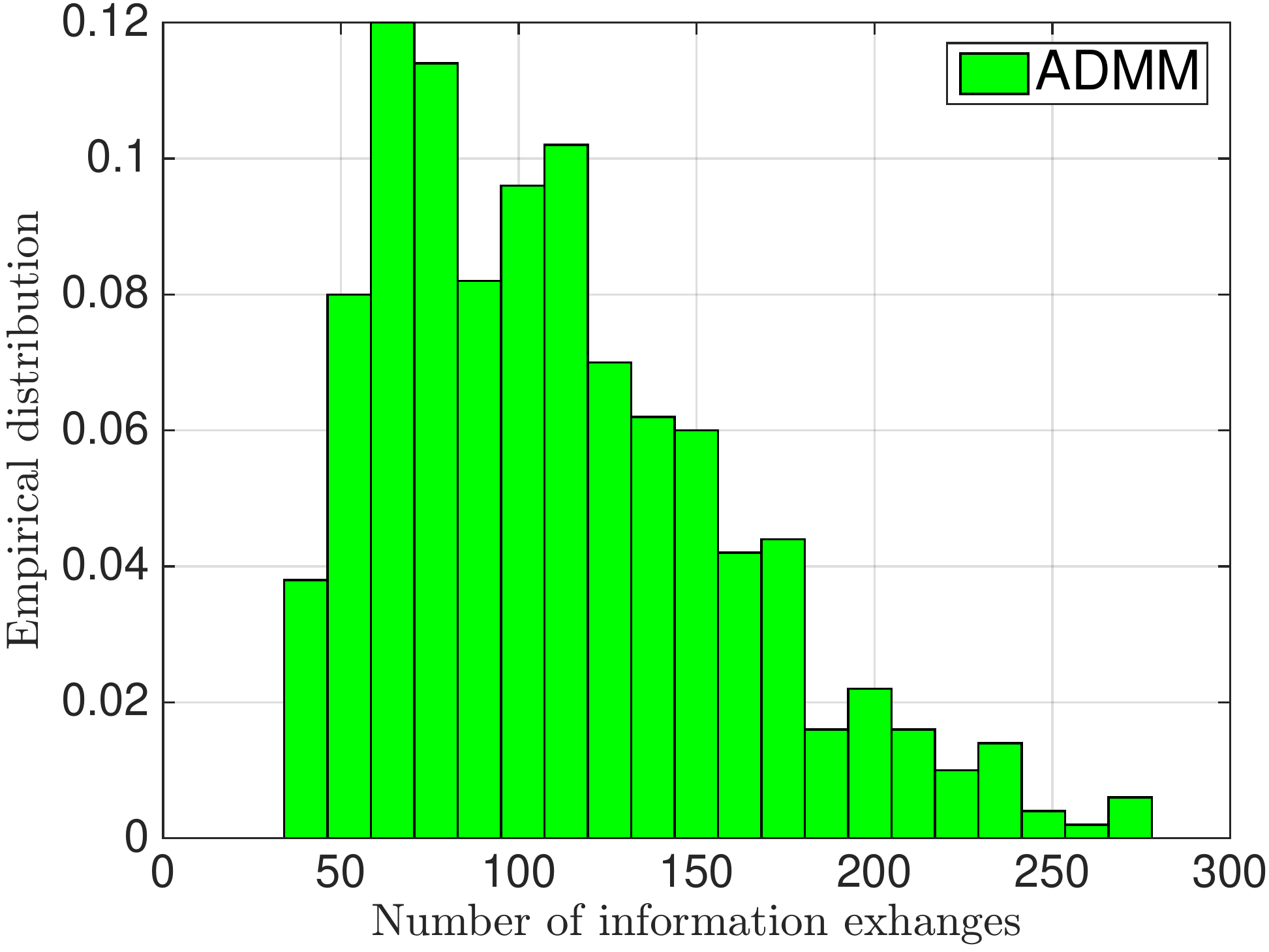}
		\caption{}\label{fig:1b}
	\end{subfigure}
	\quad
	\begin{subfigure}[t]{.3\textwidth}
		\centering
		\includegraphics[height=.13\textheight,width=\textwidth]{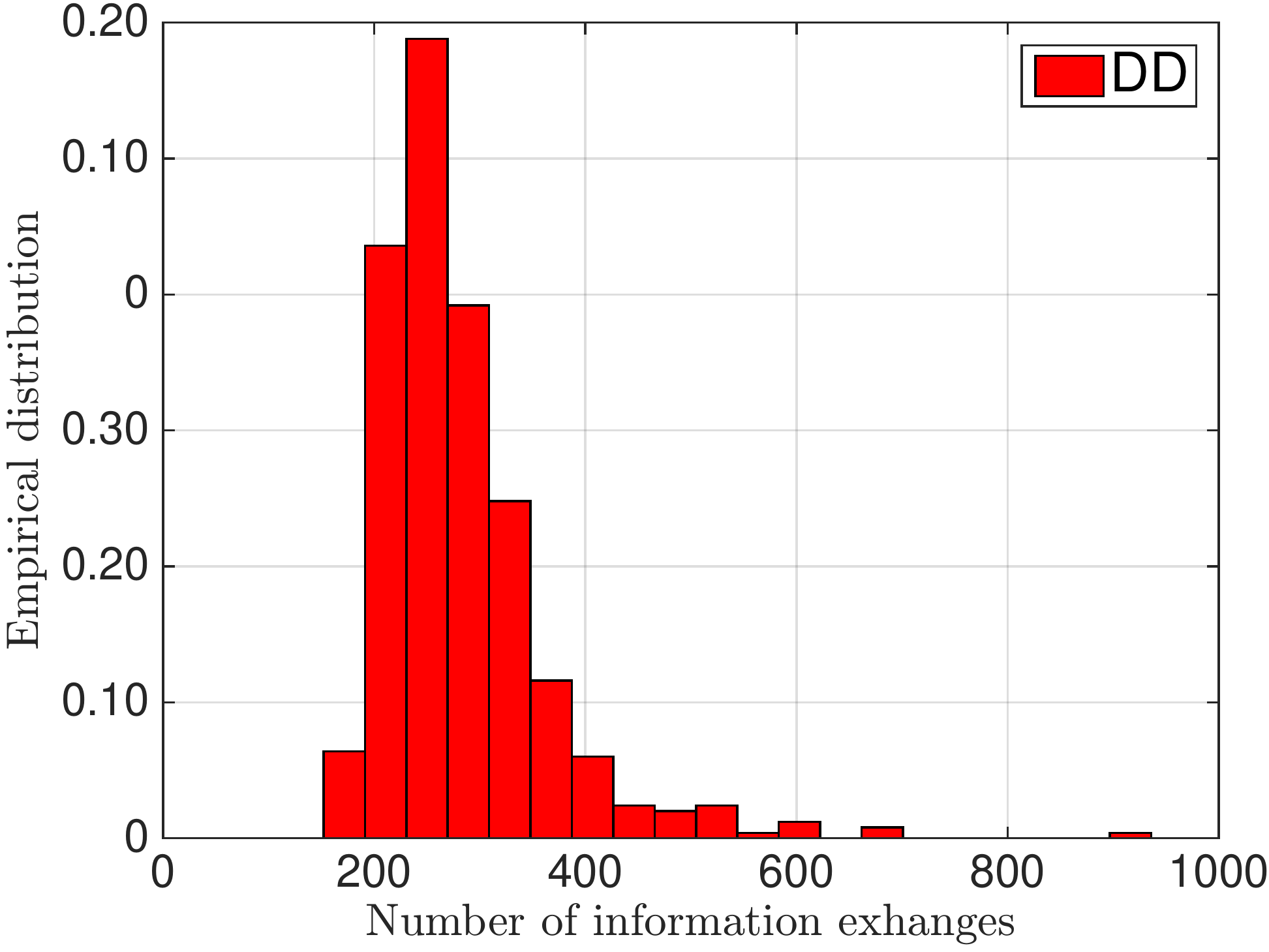}
		\caption{}\label{fig:1c}
	\end{subfigure} \\
	\begin{subfigure}[t]{.3\textwidth}
		\centering
		\includegraphics[height=.13\textheight,width=\textwidth]{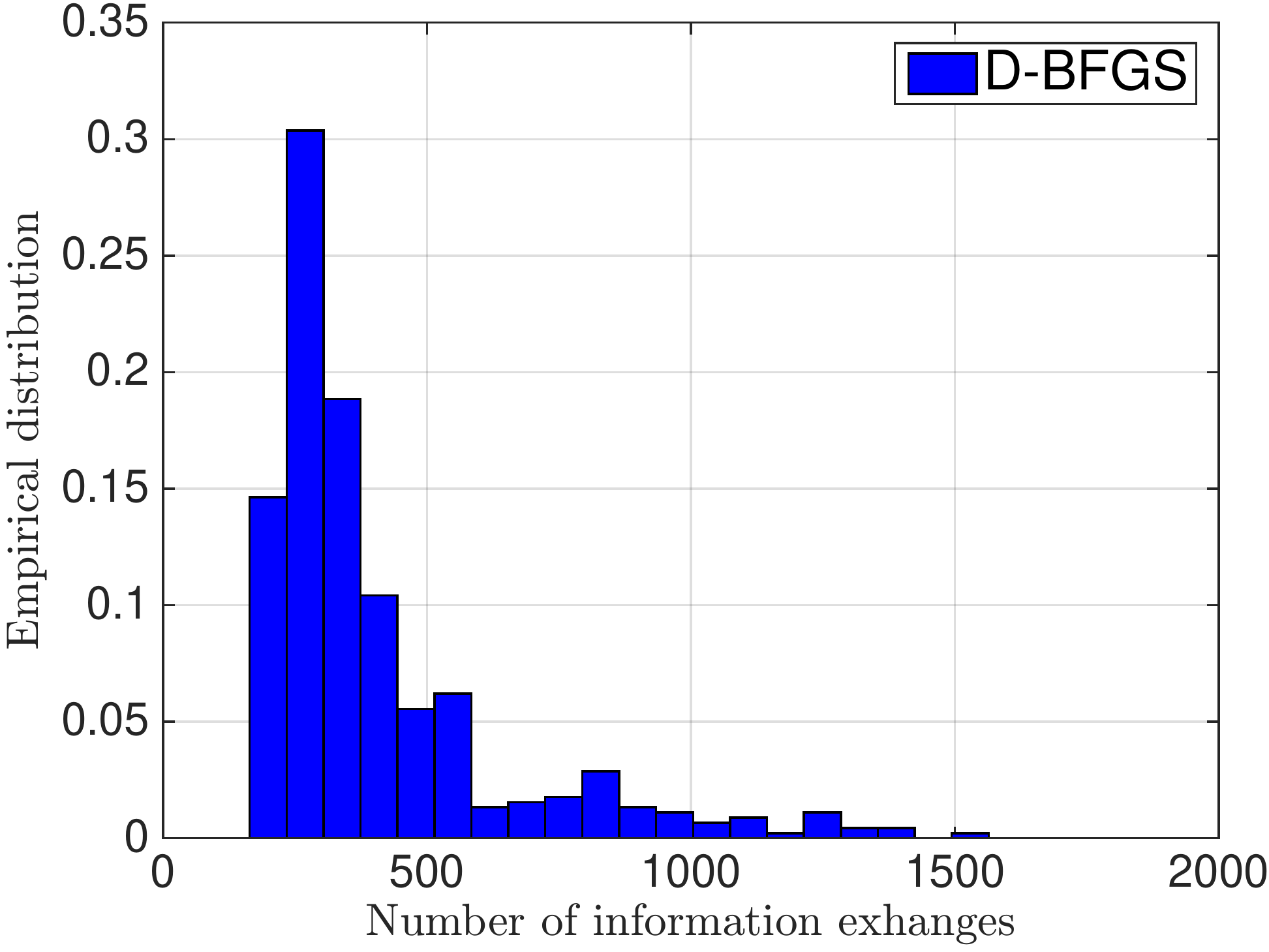}
		\caption{}\label{fig:1d}		
	\end{subfigure}
	\quad
	\begin{subfigure}[t]{.3\textwidth}
		\centering
		\includegraphics[height=.13\textheight,width=\textwidth]{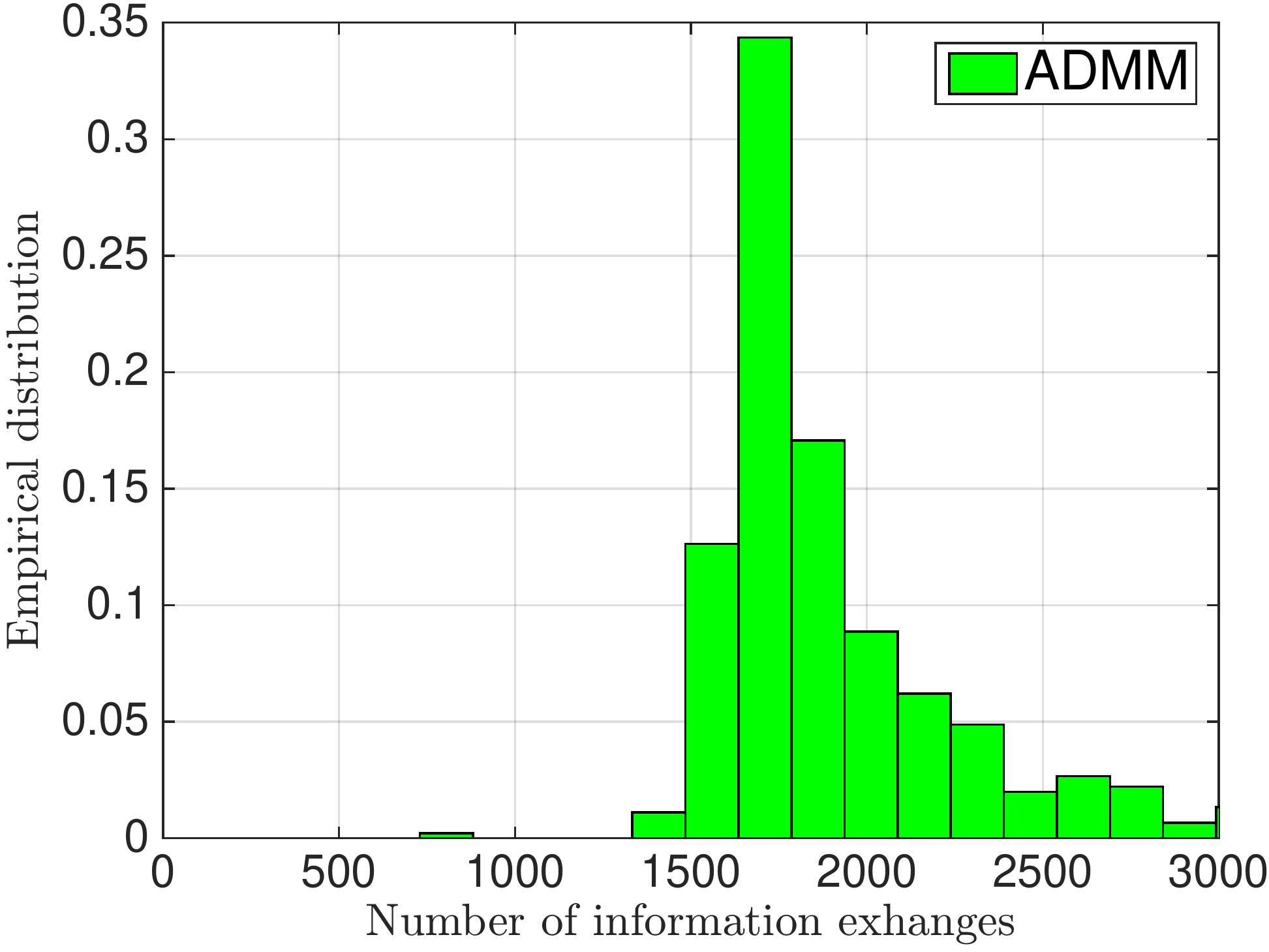}
		\caption{}\label{fig:1e}
	\end{subfigure}
	\quad
	\begin{subfigure}[t]{.3\textwidth}
		\centering
		\includegraphics[height=.13\textheight,width=\textwidth]{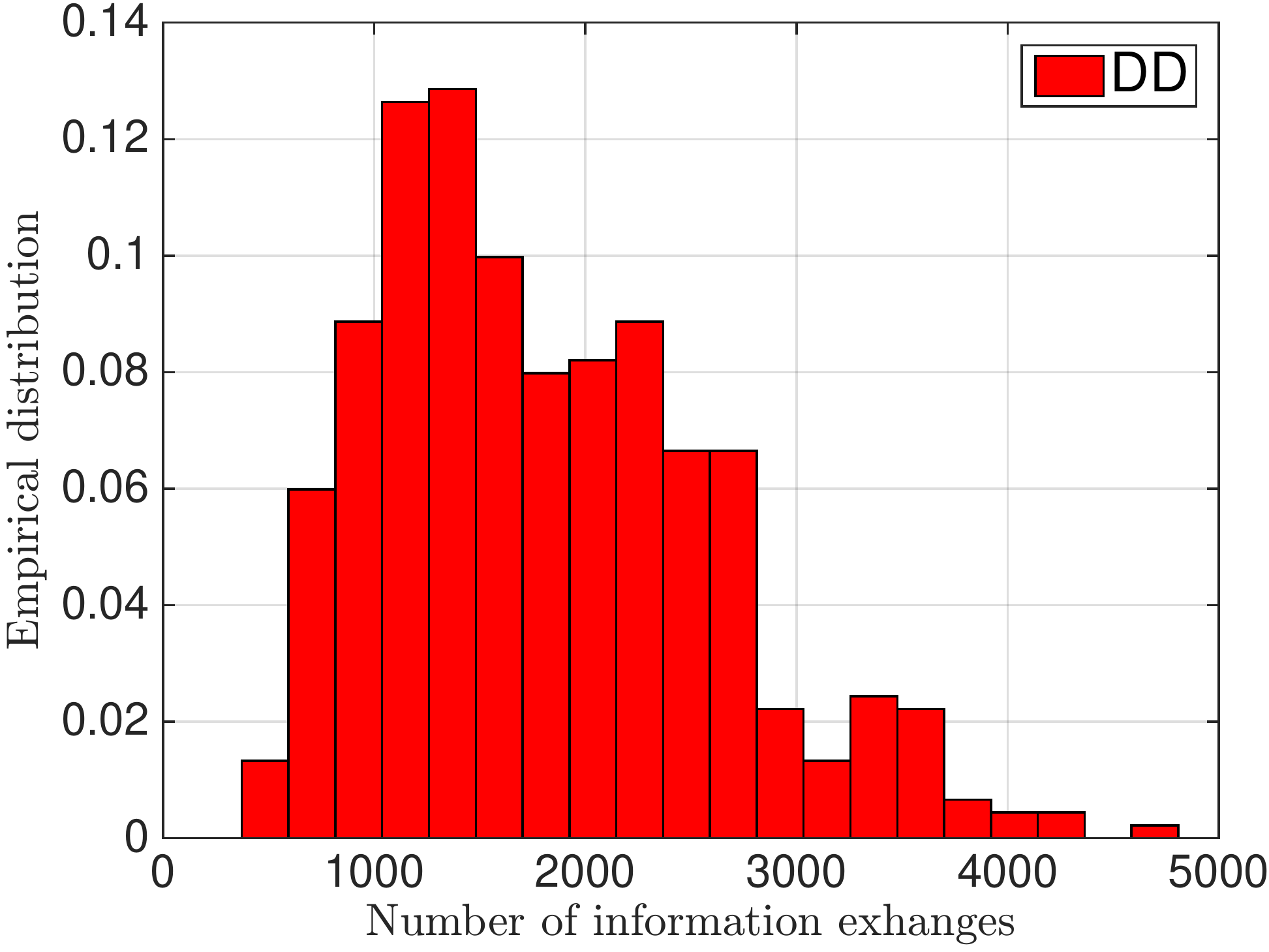}
		\caption{}\label{fig:1f}
	\end{subfigure} 	
	\caption{Histogram of number of local communications needed to converge for D-BFGS, ADMM, and DD for condition numbers of (a)-(c) $10^0$ and (d)-(f)  $10^2$. In all cases, D-BFGS provides significant improvement in convergence time over first order methods, with the larger improvement for larger condition number.}\label{figure_simulation_hist}
\end{figure*}

\begin{figure}[t]
\centering
	\begin{subfigure}[t]{.35\textwidth}
		\centering
		\includegraphics[height=.12\textheight,width=\textwidth]{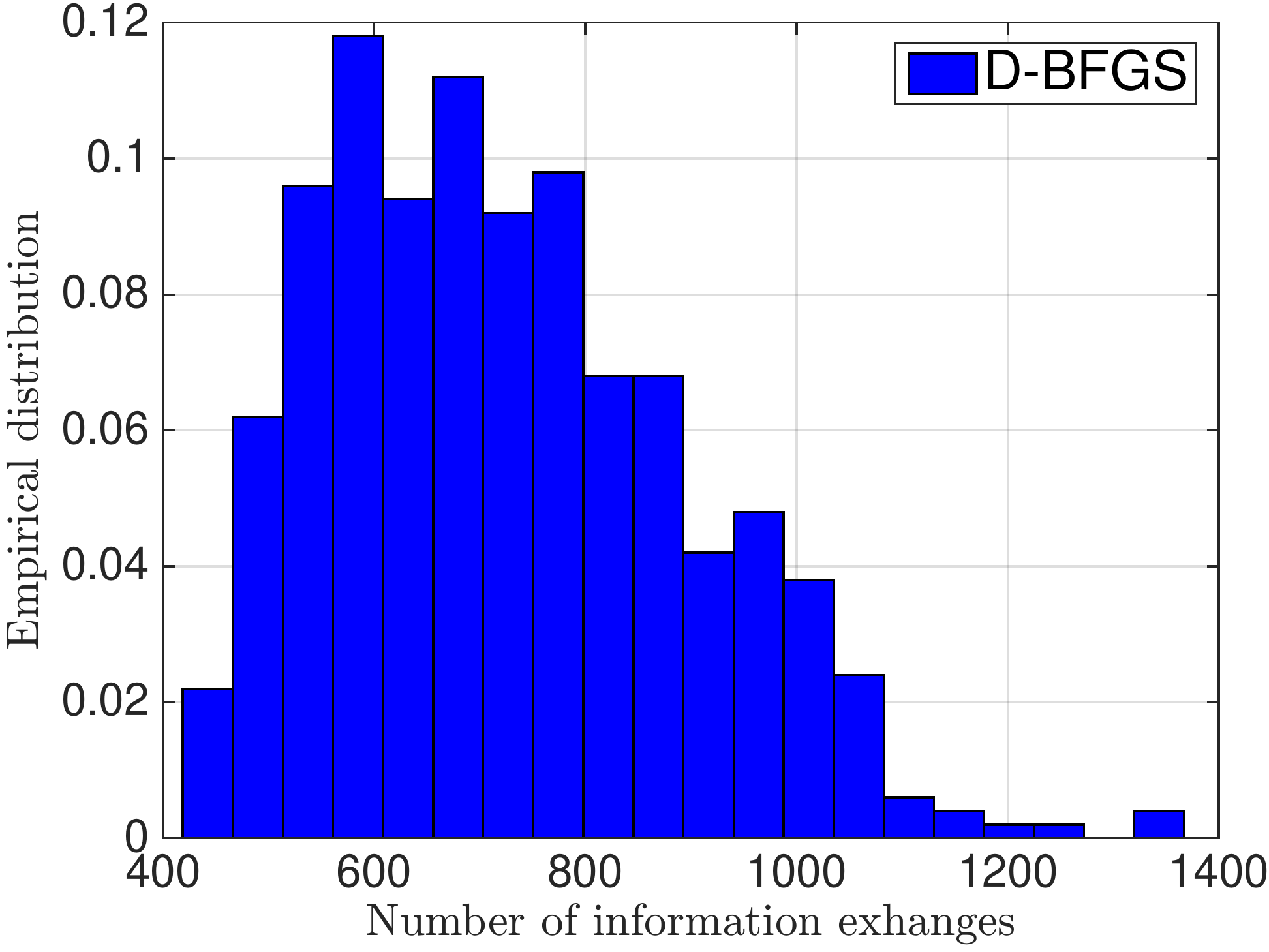}
		\caption{}\label{fig:4a}
	\end{subfigure} \\
	\begin{subfigure}[t]{.35\textwidth}
		\centering
		\includegraphics[height=.12\textheight,width=\textwidth]{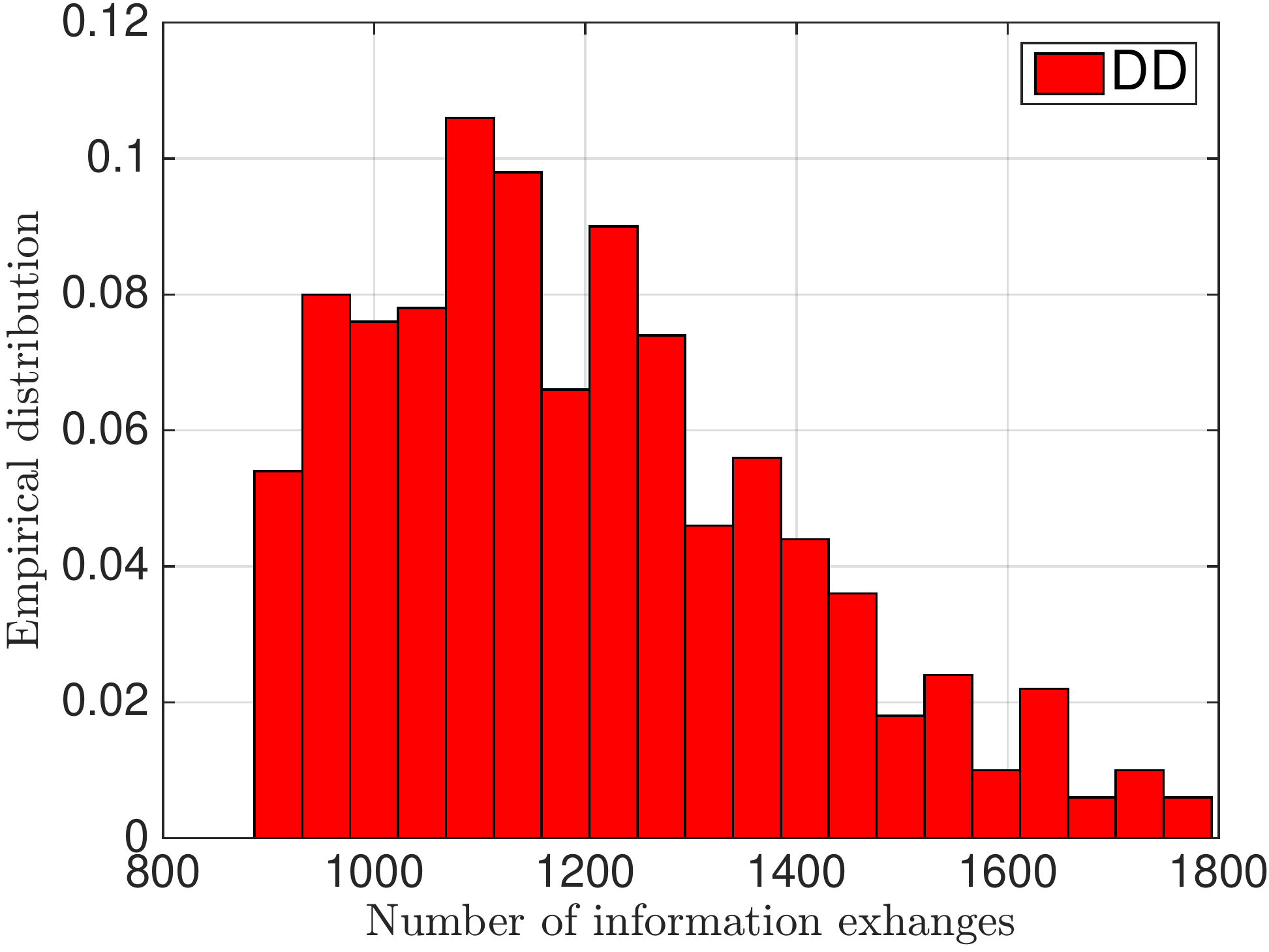}
		\caption{}\label{fig:4b}
	\end{subfigure} 	
	\caption{Histogram of number of local communications needed to converge for asynchronous (a) D-BFGS and (b) DD. D-BFGS provides significant improvement in convergence time over DD.}\label{figure_simulation_hist_async}
\end{figure}

To both ensure the local objective functions are concave and control the problem's condition number, we set the matrices $\bbA_i := \text{diag}\{ \bbA_i \}$. The first $p/2$ elements $\bba_i$
are randomly chosen from the interval $[1, 10^{-1}]$ and the last $p/2$ elements are chosen randomly from the interval $[1,10^{1}]$. The resulting $\bbA = \sum_{i=1}^n \bbA_i$ is then a positive definite with a condition number of $10^2$. For the vectors $\bbb_i$, the elements are chosen uniformly and randomly from the box $[0,1]^p$.
In our simulations we fix the variable dimension $p=4$ and the number of nodes $n=50$. The regularization parameters for D-BFGS are chosen to be $\gamma = 10^{-2}$ and $\Gamma = 10^{-3}$. In the experiments, the step size for each method is chosen to be constant and attempt is made to choose the largest step size for which the algorithms are observed to converge. For the network structure, we consider a $4$-regular cycle graph.

We simulate the performance of D-BFGS, ADMM, and DD on the dual problem in \eqref{eq_dual_problem} each using respective step sizes of $0.01$, $0.002$, and $0.002$. To view convergence with respect to the original primal problem, we look directly at how fast and how close the algorithms reach the optimal point in the original primal formulation rather than looking at the norm of dual gradient. The optimal point $\bbx^*$ is calculated for the quadratic problem in \eqref{eq_simulation_problem} using the closed form solution of a quadratic problem with linear constraints \cite{boyd} and we evaluate the normalized average error as
\begin{align}
e(t) := \frac{1}{n} \sum_{i=1}^n \frac{ \|\bbx_i(t) - \bbx^*\|^2}{ \| \bbx^*\|^2}.
\label{eq_primal_avg_error}
\end{align}

 Figure \ref{figure_simulation} shows the normalized average error $e(t)$ for all three algorithms with respect to the number of iterations. We see that D-BFGS converges substantially faster than the first order methods, reaching an average distance of $8.7 \times 10^{-5}$ by iteration 500, while ADMM and DD just reach $3.3\times10^{-2}$ and $1.8 \times10^{-1}$ respectively by iteration 500.  

It is worth noting that D-BFGS requires four local exchanges with neighbors per iteration, while ADMM and DD require only two. In Figure \ref{figure_simulation_hist} we thus present a histogram of the number of of local exchanges required for each algorithm to converge, which we define as $\delta(t) = 10^{-2}$, over the course of 1000 independent trials for both small and large condition numbers. We see that D-BFGS requires about a factor of 2 and 5 less than ADMM and DD respectively for a small condition number of $10^0$. With larger condition number, however, the improvement of D-BFGS increases to close to a factor of 7 and 8 respectively. These results showcase the advantages of decentralized quasi-Newton methods over first-order gradient descent methods. 

We additionally perform a numerical analysis on both D-BFGS and DD in the asynchronous setting for condition number $10^0$. To generate asynchronicity, we employ a simple model of aggregated Gaussian drift between nodes' local clocks. In Figure \ref{figure_simulation_hist_async}, we present histograms of number of information exchanges required to converge to an average of error of $5 \times 10^{-2}$ for both D-BFGS and DD with respective stepsizes of $0.007$ and $0.001$. While a degradation in convergence time relative to synchronous algorithms is clearly evident, D-BFGS nonetheless continues to outperform DD, requiring on average about 600 and 1200 exchanges to occur respectively.

%
\section{Conclusions} \label{sec_conclusion}
We considered the problem of decentralized consensus optimization, in which nodes sought to maximize an aggregate cost function while only being aware of a local strictly concave component. The problem was solved in the dual domain through the introduction of D-BFGS as a decentralized quasi-Newton method. In D-BFGS, a node approximates the curvature of its local cost function and its neighboring nodes to correct its descent direction. Analytical and numerical results were established showing its convergence and improvement over decentralized gradient descent methods, respectively, in synchronous and asynchronous settings.

%
\begin{appendices}

\section{Proof of Corollary \ref{cor_primal_convergence}} \label{sec_cor_primal_convergence}
The first result follows from the strong concavity of the primal function $f$, which implies that the duality gap is zero. The second result follows from the argument in \cite[Theorem 1]{beck2014gradient}, repeated here.

Consider the Lagrangian function $\ccalL(\bbx,\bblambda)$ in \eqref{eq_lagrangian}, which is strongly concave with respect to $\bbx$ with parameter $\mu$ for strongly concave $f_i$'s with parameter $\mu$. It is then the case that the following inequality holds for two points, the Lagrangian maximizer $\bbx(\bblambda)$ [cf. \eqref{eq_lagrangian_maximizers}] and the optimal primal variable $\bbx^*$,
\begin{align}
\ccalL(\bbx(\bblambda),\bblambda) - \ccalL(\bbx^*,\bblambda) \geq \frac{\mu}{2} \| \bbx(\bblambda) - \bbx^* \|^2 \label{expr09}.
\end{align}
Observe that by using the matrix $\bbA$ defined in Section \ref{sec_convergence}, we can rewrite the Lagrangian function as $\ccalL(\bbx,\bblambda) = f(\bbx) + \bblambda^T \bbA \bbx$. We can find an upper bound on $\ccalL(\bbx(\bblambda),\bblambda) - \ccalL(\bbx^*,\bblambda)$ with respect to the dual function as
\begin{align}
\ccalL(\bbx(\bblambda),\bblambda) - \ccalL(\bbx^*,\bblambda) &= f(\bbx(\bblambda)) + \bblambda^T \bbA \bbx(\bblambda) \label{expr10} \\
& \qquad - f(\bbx^*) - \bblambda^T \bbA \bbx^*. \nonumber
\end{align}
The first two terms in \eqref{expr10} are equivalent to the dual function $h(\bblambda)$, while $f(\bbx^*) = h(\bblambda^*)$ because the duality gap is zero. Additionally $\bbA \bbx^*=\bb0$ by construction so \eqref{expr10} reduces to
\begin{align}
\ccalL(\bbx(\bblambda),\bblambda) - \ccalL(\bbx^*,\bblambda) &= h(\bblambda) - h(\bblambda^*). \label{expr11}
\end{align}
We can then combine the results of \eqref{expr09}, \eqref{expr11}, and \eqref{eq_convergence_result} to obtain
\begin{align}
\frac{\mu}{2} \| \bbx(\bblambda) - \bbx^* \|^2 \leq h(\bblambda) - h(\bblambda^*) \leq o\left(\frac{1}{t}\right),
\end{align}
which provides us the with the result in \eqref{eq_primal_convergence_result}.

\end{appendices}


\bibliographystyle{IEEEtran}
\bibliography{bmc_article}

\begin{thebibliography}{10}
\providecommand{\url}[1]{#1}
\csname url@samestyle\endcsname
\providecommand{\newblock}{\relax}
\providecommand{\bibinfo}[2]{#2}
\providecommand{\BIBentrySTDinterwordspacing}{\spaceskip=0pt\relax}
\providecommand{\BIBentryALTinterwordstretchfactor}{4}
\providecommand{\BIBentryALTinterwordspacing}{\spaceskip=\fontdimen2\font plus
\BIBentryALTinterwordstretchfactor\fontdimen3\font minus
  \fontdimen4\font\relax}
\providecommand{\BIBforeignlanguage}[2]{{%
\expandafter\ifx\csname l@#1\endcsname\relax
\typeout{** WARNING: IEEEtran.bst: No hyphenation pattern has been}%
\typeout{** loaded for the language `#1'. Using the pattern for}%
\typeout{** the default language instead.}%
\else
\language=\csname l@#1\endcsname
\fi
#2}}
\providecommand{\BIBdecl}{\relax}
\BIBdecl

\bibitem{olfati2004consensus}
R.~Olfati-Saber and R.~M. Murray, ``Consensus problems in networks of agents
  with switching topology and time-delays,'' \emph{Automatic Control, IEEE
  Transactions on}, vol.~49, no.~9, pp. 1520--1533, 2004.

\bibitem{Bullo2009}
F.~Bullo, J.~Cort{\'e}s, and S.~Martinez, \emph{Distributed control of robotic
  networks: a mathematical approach to motion coordination algorithms}.\hskip
  1em plus 0.5em minus 0.4em\relax Princeton University Press, 2009.

\bibitem{Cao2013-TII}
Y.~Cao, W.~Yu, W.~Ren, and G.~Chen, ``An overview of recent progress in the
  study of distributed multi-agent coordination,'' \emph{IEEE Transactions on
  Industrial Informatics}, vol.~9, pp. 427--438, 2013.

\bibitem{LopesEtal8}
C.~G. Lopes and A.~H. Sayed, ``Diffusion least-mean squares over adaptive
  networks: Formulation and performance analysis,'' \emph{Signal Processing,
  IEEE Transactions on}, vol.~56, no.~7, pp. 3122--3136, 2008.

\bibitem{Schizas2008-1}
I.~D. Schizas, A.~Ribeiro, and G.~B. Giannakis, ``Consensus in ad hoc wsns with
  noisy links--part i: Distributed estimation of deterministic signals,''
  \emph{Signal Processing, IEEE Transactions on}, vol.~56, no.~1, pp. 350--364,
  2008.

\bibitem{KhanEtal10}
U.~A. Khan, S.~Kar, and J.~M. Moura, ``Diland: An algorithm for distributed
  sensor localization with noisy distance measurements,'' \emph{Signal
  Processing, IEEE Transactions on}, vol.~58, no.~3, pp. 1940--1947, 2010.

\bibitem{cRabbatNowak04}
M.~Rabbat and R.~Nowak, ``Distributed optimization in sensor networks,'' in
  \emph{Proceedings of the 3rd international symposium on Information
  processing in sensor networks}.\hskip 1em plus 0.5em minus 0.4em\relax ACM,
  2004, pp. 20--27.

\bibitem{bekkerman2011scaling}
R.~Bekkerman, M.~Bilenko, and J.~Langford, \emph{Scaling up machine learning:
  Parallel and distributed approaches}.\hskip 1em plus 0.5em minus 0.4em\relax
  Cambridge University Press, 2011.

\bibitem{Tsianos2012-allerton-consensus}
K.~I. Tsianos, S.~Lawlor, and M.~G. Rabbat, ``Consensus-based distributed
  optimization: Practical issues and applications in large-scale machine
  learning,'' \emph{Communication, Control, and Computing (Allerton), 2012 50th
  Annual Allerton Conference on}, pp. 1543--1550, 2012.

\bibitem{Cevher2014}
V.~Cevher, S.~Becker, and M.~Schmidt, ``Convex optimization for big data:
  Scalable, randomized, and parallel algorithms for big data analytics,''
  \emph{Signal Processing Magazine, IEEE}, vol.~31, no.~5, pp. 32--43, 2014.

\bibitem{nedic2009}
A.~Nedic and A.~Ozdaglar, ``Distributed subgradient methods for multi-agent
  optimization,'' \emph{Automatic Control, IEEE Transactions on}, vol.~54,
  no.~1, pp. 48--61, 2009.

\bibitem{nedic2010constrained}
A.~Nedi{\'c}, A.~Ozdaglar, and P.~A. Parrilo, ``Constrained consensus and
  optimization in multi-agent networks,'' \emph{Automatic Control, IEEE
  Transactions on}, vol.~55, no.~4, pp. 922--938, 2010.

\bibitem{YuanQing}
K.~Yuan, Q.~Ling, and W.~Yin, ``On the convergence of decentralized gradient
  descent,'' \emph{arXiv preprint arXiv:1310.7063}, 2013.

\bibitem{Shi2014}
W.~Shi, Q.~Ling, G.~Wu, and W.~Yin, ``Extra: An exact first-order algorithm for
  decentralized consensus optimization,'' \emph{SIAM Journal on Optimization},
  vol.~25, no.~2, pp. 944--966, 2015.

\bibitem{chatzipanagiotis2013augmented}
N.~Chatzipanagiotis, D.~Dentcheva, and M.~M. Zavlanos, ``An augmented
  lagrangian method for distributed optimization,'' \emph{Mathematical
  Programming}, pp. 1--30, 2013.

\bibitem{rockafellar1976augmented}
R.~T. Rockafellar, ``Augmented lagrangians and applications of the proximal
  point algorithm in convex programming,'' \emph{Mathematics of operations
  research}, vol.~1, no.~2, pp. 97--116, 1976.

\bibitem{jakovetic2015linear}
D.~Jakovetic, J.~M. Moura, and J.~Xavier, ``Linear convergence rate of a class
  of distributed augmented lagrangian algorithms,'' \emph{Automatic Control,
  IEEE Transactions on}, vol.~60, no.~4, pp. 922--936, 2015.

\bibitem{Jakovetic2014-1}
D.~Jakovetic, J.~Xavier, and J.~M. Moura, ``Fast distributed gradient
  methods,'' \emph{Automatic Control, IEEE Transactions on}, vol.~59, no.~5,
  pp. 1131--1146, 2014.

\bibitem{NN-part1}
A.~Mokhtari, Q.~Ling, and A.~Ribeiro, ``Network newton-part i: Algorithm and
  convergence,'' \emph{arXiv preprint arXiv:1504.06017}, 2015.

\bibitem{bajovic2015newton}
D.~Bajovic, D.~Jakovetic, N.~Krejic, and N.~K. Jerinkic, ``Newton-like method
  with diagonal correction for distributed optimization,'' \emph{arXiv preprint
  arXiv:1509.01703}, 2015.

\bibitem{mokhtari2016decentralized}
A.~Mokhtari, W.~Shi, Q.~Ling, and A.~Ribeiro, ``A decentralized second-order
  method with exact linear convergence rate for consensus optimization,''
  \emph{arXiv preprint arXiv:1602.00596}, 2016.

\bibitem{Broyden}
C.~G. Broyden, J.~E.~D. Jr., Wang, and J.~J. More, ``On the local and
  superlinear convergence of quasi-newton methods,'' \emph{IMA J. Appl. Math},
  vol.~12, no.~3, pp. 223--245, June 1973.

\bibitem{Byrd}
R.~H. Byrd, J.~Nocedal, and Y.~Yuan, ``Global convergence of a class of
  quasi-newton methods on convex problems,'' \emph{SIAM J. Numer. Anal.},
  vol.~24, no.~5, pp. 1171--1190, October 1987.

\bibitem{DingNocedal}
L.~{Dong C.} and J.~Nocedal, ``On the limited memory bfgs method for large
  scale optimization,'' \emph{Mathematical programming}, no. 45(1-3), pp.
  503--528, 1989.

\bibitem{nocedal2006numerical}
J.~Nocedal and S.~Wright, \emph{Numerical optimization}.\hskip 1em plus 0.5em
  minus 0.4em\relax Springer Science \& Business Media, 2006.

\bibitem{mokhtari2014res}
A.~Mokhtari and A.~Ribeiro, ``Res: Regularized stochastic bfgs algorithm,''
  \emph{Signal Processing, IEEE Transactions on}, vol.~62, no.~23, pp.
  6089--6104, 2014.

\bibitem{bertsekas1989parallel}
D.~P. Bertsekas and J.~N. Tsitsiklis, \emph{Parallel and distributed
  computation: numerical methods}.\hskip 1em plus 0.5em minus 0.4em\relax
  Prentice-Hall, Inc., 1989.

\bibitem{bertsekas1999nonlinear}
D.~P. Bertsekas, ``Nonlinear programming,'' 1999.

\bibitem{DBFGS2}
M.~Eisen, A.~Mokhtari, and A.~Ribeiro, ``Decentralized quasi-newton methods,''
  2016, available at {\footnotesize
  {\href{http://www.seas.upenn.edu/~maeisen/wiki/dbfgs.pdf}{http://www.seas.upenn.edu/$\sim$maeisen/wiki/dbfgs.pdf}}}.

\bibitem{boyd}
S.~Boyd and L.~Vandenberghe, \emph{Convex optimization}.\hskip 1em plus 0.5em
  minus 0.4em\relax Cambridge university press, 2004.

\bibitem{beck2014gradient}
A.~Beck, A.~Nedic, A.~Ozdaglar, and M.~Teboulle, ``An gradient method for
  network resource allocation problems,'' \emph{Control of Network Systems,
  IEEE Transactions on}, vol.~1, no.~1, pp. 64--73, 2014.

\end{thebibliography}

\end{document}